# Derivative Polynomials for tanh, tan, sech and sec in Explicit Form

Khristo N. Boyadzhiev
Ohio Northern University
Department of Mathematics
Ada, OH 45810

k-boyadzhiev@onu.edu

**Abstract.**
The derivative polynomials for the hyperbolic and trigonometric tangent, cotangent and secant are found in explicit form, where the coefficients are given in terms of the Stirling numbers of the second kind. As application we evaluate some integrals and also give the reflection formula for the polygamma function in explicit form.

## 0. Introduction.

The derivative polynomials for tangent and secant, correspondingly $P_m$ and $Q_m$, are defined by the equations ($m = 0, 1, \ldots$):

$$(\frac{d}{d\theta})^m \tan\theta = P_m(\tan\theta), \qquad (0.1)$$

$$(\frac{d}{d\theta})^m \sec\theta = \sec\theta \, Q_m(\tan\theta)$$

(see Hoffman [5]). Thus $P_m(z)$ is a polynomial of degree $m+1$ and $Q_m(z)$ is a polynomial of degree $m$. The derivative polynomials have generating functions

$$\tan(\theta + x) = \sum_{m=0}^{\infty} P_m(\tan\theta) \frac{x^m}{m!},$$

$$\frac{\sec(\theta+x)}{\sec(\theta)} = \sum_{m=0}^{\infty} Q_m(\tan\theta) \frac{x^m}{m!},$$

which result from Taylor's expansions centered at $x=0$ of the function $f(x) = \tan(\theta+x)$ and $g(x) = \sec(\theta+x)$ for a fixed $\theta$ (cf. [4, p.287], [5], [8]). Recurrence relations are known for both $P_m$ and $Q_m$, see [5]. Hoffman presented a number of properties of these polynomials in [5] and [6], but did not compute them explicitly. Another study on the derivative polynomials can be found in Schwatt [13].

In this note we provide an explicit form of the derivative polynomials for tanh, tan, sech and sec, where the coefficients are written in terms of Stirling numbers and $P_m(0), Q_m(0)$ are expressed in terms of Bernoulli and Euler numbers, respectively. We also give some applications, complementing the interesting applications in [5].

It is well-known that the Eisenstein series

$$e_r(z) = \frac{1}{z^r} + \sum_{k=1}^{+\infty} \left( \frac{1}{(z+k)^r} + \frac{1}{(z-k)^r} \right)$$

($r$ - a positive integer) can be written as a polynomial of $e_1(z)$, see [10, Chapter 11]. In Section 4 we give this polynomial in explicit form (equation (4.3)). In that section we consider also the Polygamma functions

$$\psi_n(z) := \left(\frac{d}{dz}\right)^{n+1} \log\Gamma(z),$$

and provide an explicit form for their reflection formula:

$$\psi_n(z) - (-1)^n \psi_n(1-z) = -\pi \left(\frac{d}{dz}\right)^n \cot\pi z. \tag{0.2}$$

Finally, in Section 5 the derivative polynomials are used to evaluate some integrals.

# 1. Prerequisites.

Let $\{{n \atop k}\}$ be the Stirling numbers of the second kind; see e.g., [4, Ch 6]. Consider the geometric polynomials

$$\omega_n(x) = \sum_{k=0}^{n} \{{n \atop k}\} k! \, x^k, \tag{1.1}$$

which serve as derivative polynomials for the geometric series

$$\frac{1}{1-x} = \sum_{k=0}^{\infty} x^k, \quad |x| < 1,$$

according to the derivative operator $xD$ ($D = d/dx$), that is,

$$(xD)^m \left\{\frac{1}{1-x}\right\} = \frac{1}{1-x} \omega_m\left(\frac{x}{1-x}\right), \quad m = 0, 1, 2, \ldots, \quad x \neq 1 \tag{1.2a}$$

or

$$\sum_{k=0}^{\infty} k^m x^k = \frac{1}{1-x} \omega_m\left(\frac{x}{1-x}\right), \text{ when } |x| < 1, \tag{1.2b}$$

see [2], [3, (4.10)].

We define now a more general class of polynomials. For any $a, b \in \mathbb{C}$ set

$$p_m(z; a, b) := \sum_{k=0}^{m} \binom{m}{k} a^{m-k} b^k \omega_k(z), \tag{1.3a}$$

or

$$p_m(z; a, b) = \sum_{j=0}^{m} \left[ j! \sum_{k=j}^{m} \binom{m}{k} \{{k \atop j}\} a^{m-k} b^k \right] z^j. \tag{1.3b}$$

Note that when $a = 0$:

$$p_m(z; 0, b) = b^m \omega_m(z). \tag{1.4}$$

The importance of these polynomials comes from the following lemma.

**Lemma 1.** For every $a, b \in \mathbb{C}$ with $x^b \neq 1$ and every $m = 0, 1, \ldots$, we have:

$$(xD)^m \frac{x^a}{1-x^b} = \frac{x^a}{1-x^b} p_m\left(\frac{x^b}{1-x^b}; a, b\right), \tag{1.5a}$$

$$(xD)^m \frac{x^a}{1+x^b} = \frac{x^a}{1+x^b} P_m(\frac{-x^b}{1+x^b};a,b). \tag{1.5b}$$

*Proof.* We prove (1.5a). Let $|x^b|<1$. Then

$$(xD)^m \frac{x^a}{1-x^b} = (xD)^m \sum_{n=0}^{\infty} x^{a+bn} = \sum_{n=0}^{\infty} (a+bn)^m x^{a+bn}$$

$$= x^a \sum_{n=0}^{\infty} \sum_{k=0}^{m} \binom{m}{k} a^{m-k} b^k n^k x^{bn} = x^a \sum_{k=0}^{m} \binom{m}{k} a^{m-k} b^k \{\sum_{n=0}^{\infty} n^k (x^b)^n\}$$

$$= x^a \sum_{k=0}^{m} \binom{m}{k} a^{m-k} b^k \{\frac{1}{1-x^b} \omega_k(\frac{x^b}{1-x^b})\} = \frac{x^a}{1-x^b} P_m(\frac{x^b}{1-x^b};a,b).$$

Now (1.5a) follows for every $x^b \neq 1$ by analytic continuation. (1.5b) is proved in the same manner.

## 2. The derivative polynomials for *coth*, *tanh*, *tan* and *cot*.

Define the derivative polynomial $C_m(z)$ of degree $m+1$ by

$$(\frac{d}{d\theta})^m \coth\theta = C_m(\coth\theta). \tag{2.1}$$

Therefore, $C_0(z) = 1$, $C_1(z) = 1-z^2$, $C_2(z) = 2z^3 - 2z$, etc.

**Proposition 1.** For $m \geq 1$ one has

$$C_m(z) = (-2)^m (z+1) \omega_m(\frac{z-1}{2}), \tag{2.2}$$

where $\omega_m$ is the geometric polynomial (1.1).

*Proof.* We have

$$\coth\theta = \frac{e^\theta + e^{-\theta}}{e^\theta - e^{-\theta}} = \frac{2}{1-e^{-2\theta}} - 1.$$

Let $x = e^\theta$. Then $d/d\theta = xD$ and from Lemma 1,

$$(\frac{d}{d\theta})^m \coth \theta = 2(xD)^m \frac{1}{1-x^{-2}} = 2(-2)^m \frac{1}{1-x^{-2}} \omega_m(\frac{x^{-2}}{1-x^{-2}}).$$

Set $z = \coth \theta$. Then

$$\frac{1}{1-e^{-2\theta}} = \frac{e^\theta}{e^\theta - e^{-\theta}} = \frac{\cosh \theta + \sinh \theta}{2 \sinh \theta} = \frac{1}{2}(z+1),$$

and,

$$\frac{e^{-2\theta}}{1-e^{-2\theta}} = \frac{e^{-\theta}}{e^\theta - e^{-\theta}} = \frac{\cosh \theta - \sinh \theta}{2 \sinh \theta} = \frac{1}{2}(z-1),$$

so we arrive at (2.2).

The higher derivatives of **tanh** are formed in the same pattern as those of **coth**, as both functions satisfy the same differential equation

$$y' = 1 - y^2.$$

Therefore, $C_m$ are also the derivative polynomials of **tanh**,

$$(\frac{d}{d\theta})^m \tanh \theta = C_m(\tanh \theta). \tag{2.3}$$

Explicitly, for $m \geq 1$

$$C_m(z) = (-2)^m (z+1) \sum_{k=0}^{m} \frac{k!}{2^k} \{{}^m_k\} (z-1)^k. \tag{2.4}$$

Thus $C_m(-1) = 0$, $C_m(1) = 2(-2)^m \{{}^m_0\} = 0$ ($m > 0$) and

$$C_m(0) = (-2)^m \omega_m(\frac{-1}{2}) = (-2)^m \sum_{k=0}^{m} (-1)^k \{{}^m_k\} \frac{k!}{2^k}. \tag{2.5}$$

The numbers $C_m(0)/m!$ are the coefficients in the Maclaurin expansion of $\tanh \theta$. In view of [4, p. 317, problem 6.76] one has

$$\sum_{k=0}^{m} (-1)^k \{{}^m_k\} \frac{k!}{2^k} = \frac{2}{m+1}(1-2^{m+1})B_{m+1},$$

where $B_k$ are the Bernoulli numbers, and then from (2.5)

$$C_m(0) = \frac{(-1)^m}{m+1} 2^{m+1}(1-2^{m+1})B_{m+1}. \tag{2.6}$$

Note that for even indices $m = 2k$ the numbers $C_{2k}(0)$ are zero as the hyperbolic tangent is an odd function. When $m = 2k-1$:

$$C_{2k-1}(0) = \frac{1}{2k} 2^{2k}(2^{2k}-1)B_{2k}. \tag{2.7}$$

**Proposition 2**. For $m \geq 1$ the derivative polynomials for the tangent are:

$$P_m(z) = i^{m+1} 2^m (1-iz) \omega_m(\frac{1+iz}{-2}). \tag{2.8}$$

This result follows, in the same way as in the proof of Proposition 1, from the representation

$$\tan\theta = \frac{2i}{1+e^{2i\theta}} - i.$$

We can obtain $P_m(z)$ also from Proposition 1 by using the relation

$$\tan(\theta) = -i\tanh(i\theta).$$

Thus

$$P_m(\tan\theta) = (\frac{d}{d\theta})^m \tan\theta = -i^{m+1} C_m(\tanh(i\theta)) = -i^{m+1} C_m(i\tan\theta), \tag{2.9}$$

and from (2.2)

$$P_m(z) = -i^{m+1}(-2)^m (iz+1) \omega_m(\frac{iz-1}{2}), \tag{2.10}$$

which is the same as (2.8), as $\tan\theta$ is an odd function and from (0.1)

$$P_m(-z) = (-1)^{m+1} P_m(z).$$

In view of (1.1),

$$P_m(z) = -i^{m+1}(-2)^m (iz+1) \sum_{k=0}^{m} \frac{k!}{2^k} \{{}^m_k\} (iz-1)^k, \qquad (2.11)$$

and in particular,

$$P_m(0) = -i^{m+1}(-2)^m \sum_{k=0}^{m} (-1)^k \{{}^m_k\} \frac{k!}{2^k} \qquad (2.12)$$

(obtained independently by Knopf [7]). Here $P_{2k}(0) = 0$ and the numbers

$$P_{2k-1}(0) = (-1)^{k+1} C_{2k-1}(0) \qquad (2.13)$$

(see (2.9)) are known as the *tangent numbers*. In view of (2.7) we have the well-known representation of these numbers ([4, p.287], [11, p. 154], [10, p.221]):

$$P_{2k-1}(0) = \frac{(-1)^{k+1}}{2k} 2^{2k}(2^{2k}-1) B_{2k}. \qquad (2.14)$$

Finally, the derivative polynomials for the cotangent are

$$-P_m(-z) = (-1)^m P_m(z). \qquad (2.15)$$

This follows immediately from the equation

$$\cot\theta = i\coth(i\theta).$$

### 3. The derivative polynomials for *sech* and *sec*.

Define the polynomials $S_m(z)$ of degree $m = 0, 1, \ldots$, by

$$\left(\frac{d}{d\theta}\right)^m \operatorname{sech}\theta = \operatorname{sech}\theta\, S_m(\tanh\theta). \qquad (3.1)$$

One easily computes $S_0(z) = 1$, $S_1(z) = -z$, $S_2(z) = 2z^2 - 1$, etc.

**Proposition 3.** For $m = 0, 1, \ldots$,

$$S_m(z) = p_m(\frac{1+z}{-2}; 1, 2) \tag{3.2}$$

$$= \sum_{k=0}^{m} \binom{m}{k} 2^k \omega_k(\frac{z+1}{-2}), \text{ or}$$

$$S_m(z) = \sum_{j=0}^{m} [(-1)^j j! \sum_{k=j}^{m} \binom{m}{k}\{{}^k_j\} 2^{k-j}](z+1)^j$$

*Proof.* Starting from the representation

$$\text{sech}\,\theta = \frac{2e^\theta}{1+e^{2\theta}},$$

we apply Lemma 1 with $x = e^\theta$ to get

$$(\frac{d}{d\theta})^m \text{sech}\,\theta = 2\,(xD)^m \frac{x}{1+x^2} = \frac{2x}{1+x^2} p_m(\frac{-x^2}{1+x^2}; 1, 2)$$

$$= \text{sech}\,\theta\, p_m(\frac{-e^{2\theta}}{1+e^{2\theta}}; 1, 2),$$

which is (3.2), because

$$\frac{-e^{2\theta}}{1+e^{2\theta}} = \frac{-1}{2}(1 + \tanh\theta).$$

**Remark**. The numbers $E_m = S_m(0)$ are the *Euler numbers*, defined by

$$\text{sech}(z) = \sum_{n=0}^{\infty} E_m \frac{z^m}{n!}$$

(see [1], [11, p.153]). From (3.2) and (2.6):

$$E_m = \sum_{k=0}^{m} \binom{m}{k} 2^k \omega_k(\frac{-1}{2}) = \sum_{k=0}^{m} \binom{m}{k}(-1)^k C_k(0), \text{ or}$$

$$E_m = \sum_{k=0}^{m} \binom{m}{k} \frac{2^{k+1}}{k+1} (1-2^{k+1})B_{k+1}, \text{ or}$$

$$E_n = \frac{1}{m+1} \sum_{k=1}^{m+1} \binom{m+1}{k} 2^k (1-2^k) B_k \qquad (3.3)$$

cf. [11, p.152].

**Remark.** $S_m(z)$ is also the derivative polynomial for the hyperbolic cosecant, i.e.,

$$(\frac{d}{d\theta})^m \operatorname{csch} \theta = \operatorname{csch} \theta \, S_m(\coth \theta). \qquad (3.4)$$

This can be seen easily by the same method using the representation

$$\operatorname{csch} \theta = \frac{-2e^\theta}{1 - e^{2\theta}},$$

or by looking at the formation of the derivatives. It is interesting that if one uses the alternative representation

$$\operatorname{csch} \theta = \frac{2e^{-\theta}}{1 - e^{-2\theta}} \qquad (3.5)$$

for the above proof, the result is an equivalent form of $S_m(z)$.

We turn now to the derivative polynomials $Q_m$ for the secant,

$$(\frac{d}{d\theta})^m \sec \theta = \sec \theta \, Q_m(\tan \theta).$$

**Proposition 4.** For $m = 1, 2, \ldots$ one has:

$$Q_m(z) = i^m S_m(iz). \qquad (3.6)$$

The proof follows from (3.1) and the equation

$$\sec \theta = \operatorname{sech}(i\theta).$$

It is good to mention here also a convenient differentiation rule for the cosecant which follows from the representation

$$\frac{1}{\sin\theta} = \frac{1}{2}(\tan\frac{\theta}{2} + \cot\frac{\theta}{2}), \qquad (3.7)$$

$$(\frac{d}{d\theta})^m \frac{1}{\sin\theta} = \frac{1}{2^{m+1}}[P_m(\tan\frac{\theta}{2}) + (-1)^m P_m(\cot\frac{\theta}{2})]. \qquad (3.8)$$

**4. The Hilbert-Eisenstein series and the reflection formula for the Polygamma function.**

Consider the series

$$e_r(z) = \frac{1}{z^r} + \sum_{k=1}^{+\infty}(\frac{1}{(z+k)^r} + \frac{1}{(z-k)^r}),$$

when $r > 0$, or

$$e_r(z) = \sum_{k=-\infty}^{+\infty}\frac{1}{(z+k)^r}$$

when $r > 1$. Its theory is nicely presented in Chapter 11 of Remmert's book [10], where also a biographical note on Ferdinand Gotthold Max Eisenstein (1823-1852) can be found. We shall evaluate this series using the tangent derivative polynomials. It is known that [10, p.326]:

$$e_1(z) = \pi\cot\pi z. \qquad (4.1)$$

Differentiating this equation $r - 1$ times we find according to (2.15)

$$(-1)^{r-1}(r-1)!\, e_r(z) = \pi(\frac{d}{dz})^{r-1}\cot\pi z = -\pi^r P_{r-1}(-\cot\pi z), \qquad (4.2)$$

and therefore, from (2.11),

$$e_r(z) = \frac{-i^r\pi^r}{(r-1)!}\csc^2\pi z \sum_{k=1}^{r-1}\{_k^{r-1}\}2^{r-k-1}k!(i\cot\pi z - 1)^{k-1}.$$

We can also write $e_m$ explicitly as a polynomial of $e_1$. From (4.1) and (4.2),

$$e_r = \frac{(-1)^r\pi^r}{(r-1)!}P_{r-1}(\frac{-1}{\pi}e_1). \qquad (4.3)$$

Next we look at the digamma function

$$\psi(z) = \frac{d}{dz} \log \Gamma(z) \tag{4.4}$$

where $\Gamma(z)$ is Euler's Gamma function. The Polygamma functions $\psi_n$ are defined as the higher derivatives of $\psi$, see [1, 6.4], [12]. Thus $\psi_0(z) = \psi(z)$ and

$$\psi_n(z) := \psi^{(n)}(z) = \left(\frac{d}{dz}\right)^{n+1} \log \Gamma(z)$$

($n = 1, 2, \ldots$). The Digamma function satisfies the equation [1], [12, p.14]:

$$\psi(z) - \psi(1-z) = -\pi \cot \pi z,$$

(which reveals, in particular, that $\psi(z) - \psi(1-z) = -e_1(z)$). Differentiating this equation $n$ times, we find the exact form of the *reflection formula* for the Polygamma function [1, p.260], [12, p.22], [14]:

$$\psi_n(z) - (-1)^n \psi_n(1-z) = \pi^{n+1} P_n(-\cot \pi z) \tag{4.5}$$

$$= -i^{n+1} \pi^{n+1} \csc^2 \pi z \sum_{k=1}^{n} \left\{{n \atop k}\right\} 2^{n-k} k! (i \cot \pi z - 1)^{k-1}.$$

## 5. Evaluation of integrals.

Hoffman presented in [5, theorems 4.1 and 4.2] the following evaluations

$$\int_{-\infty}^{+\infty} \frac{x^n e^{ax}}{e^x + 1} dx = \pi^{n+1} \csc a\pi \, Q_n(-\cot a\pi), \; n \geq 0,$$

$$\int_{-\infty}^{+\infty} \frac{x^n e^{ax}}{e^x - 1} dx = \pi^{n+1} P_n(-\cot a\pi), \; n \geq 1,$$

where $0 < a < 1$. We want to add some more examples.

**Example 1.** Consider the integral [9, 2.5.46.5, p. 468]

$$\int_0^\infty \frac{\cos ax}{\cosh x} dx = \frac{\pi}{2} \operatorname{sech} \frac{\pi a}{2}, \tag{5.1}$$

for all real $a$. Differentiating this equation $n$ times with respect to $a$ we find

$$\int_0^\infty \frac{x^n}{\cosh x} \cos(ax + \frac{\pi n}{2}) dx = (\frac{\pi}{2})^{n+1} \operatorname{sech} \frac{\pi a}{2} S_n(\tanh \frac{\pi a}{2}). \tag{5.2}$$

**Example 2.** Next we consider a similar integral [9, 2.5.46.4, p. 467]:

$$\int_0^\infty \frac{\sin ax}{\cosh x} dx = -\frac{\pi}{2} \tanh \frac{\pi a}{2} + \frac{i}{2} [\psi(\frac{1-ia}{4}) - \psi(\frac{1+ia}{4})], \tag{5.3}$$

where $\psi(z)$ is the Digamma function (4.4). This leads to

$$\int_0^\infty \frac{x^n}{\cosh x} \sin(ax + \frac{\pi n}{2}) dx = -(\frac{\pi}{2})^{n+1} C_n(\tanh \frac{\pi a}{2}) \tag{5.4}$$

$$+ \frac{i}{2^{2n+1}} [(-i)^n \psi^{(n)}(\frac{1-ia}{4}) - i^n \psi^{(n)}(\frac{1+ia}{4})].$$

In the same way one obtains explicit evaluations of the integrals (2) and (5) in section 2.5.47, p. 469 in [9].

The two integrals (5.2) and (5.4) can be united in one to give a more compact result.

**Proposition 5.** For every $n = 0, 1, \ldots$ and every real $a$,

$$\int_0^\infty \frac{x^n}{\cosh x} e^{iax} dx = \tag{5.5}$$

$$(-1)^n i^n (\frac{\pi}{2})^{n+1} [\operatorname{sech} \frac{\pi a}{2} S_n(\tanh \frac{\pi a}{2}) - i C_n(\tanh \frac{\pi a}{2})]$$

$$+ \frac{1}{2^{2n+1}}[\psi^{(n)}(\frac{1+ia}{4}) - (-1)^n \psi^{(n)}(\frac{1-ia}{4})].$$

This formula obviously extends to complex $|a| < 1$ (cf. 2.5.47(4), p.469 in [9]).

At the end, the author wishes to express his gratitude to the referee for a number of helpful remarks.

### References.